\documentclass[12pt]{amsart}
\usepackage{amsmath,amssymb,xy,verbatim,amscd}

\swapnumbers
\theoremstyle{plain}
\newtheorem{theorem}{Theorem}[section]

\newtheorem{proposition}[theorem]{Proposition}

\theoremstyle{definition}
\newtheorem{definition}[theorem]{Definition}

\newtheorem{conjecture}[theorem]{Conjecture}
\newtheorem{subsect}[theorem]{}

\numberwithin{equation}{theorem}

\allowdisplaybreaks[1]

\newcommand{\ds}{\displaystyle}

\DeclareMathOperator{\Hom}{Hom}

\newcommand\C{\mathbb{C}}

\newcommand\frh{\mathfrak{h}}
\newcommand\frg{\mathfrak{g}}
\newcommand\frb{\mathfrak{b}}

\begin{document}

\title{A Conjectural Presentation of Fusion Algebras}
\author{Arzu Boysal and Shrawan Kumar}
\date{}
\maketitle \pagestyle{myheadings}
{\it Dedicated to Prof. Masaki Kashiwara on his sixtieth birthday}

\markboth{}{}  \maketitle

\section{Introduction}
Let $G$ be a connected, simply-connected, simple  algebraic
group over $\C$. We fix a Borel subgroup $B$ of $G$
 and a maximal torus $T\subset B$. We denote their Lie algebras by $\frg,
 \mathfrak b, \mathfrak h$ respectively.  Let $P_+\subset \mathfrak h^*$ be the set of
 dominant integral weights. For any
 $\lambda \in P_+,$ let $V(\lambda)$ be the finite dimensional
 irreducible $\frg$-module
 with highest weight $\lambda$.  We fix a positive
integer $\ell$ and let $\mathcal{R}_{\ell}(\frg)$ be the free
$\mathbb{Z}$-module with basis $\{ V({\lambda}): \lambda \in
P_{\ell} \}$, where
\[P_\ell:=\{\lambda\in P_+: \lambda(\theta^\vee)\leq \ell\},\]
 $\theta$ is the highest root of $\frg$ and $\theta^\vee$ is the associated coroot.
There is a product structure on $\mathcal{R}_{\ell}(\frg)$, called
the {\it fusion product} (cf. Section 3), making it a commutative associative
(unital) ring.  In this paper, we consider its
complexification $\mathcal{R}_{\ell}^\C(\frg)$, called the {\it
fusion algebra}, which is a finite dimensional (commutative and associative) algebra
without nilpotents.

Let $\mathcal{R}(\frg)$ be the Grothendieck ring of finite
dimensional representations of $\frg$ and let $\mathcal{R}^\C(\frg)$ be its complexification.
As in subsection \ref{3.4}, there is a surjective ring homomorphism
$\beta: \mathcal{R}(\frg) \to \mathcal{R}_\ell(\frg)$. Let $\beta^\C:
\mathcal{R}^\C(\frg) \to \mathcal{R}_\ell^\C(\frg)$ be its complexification and let
 $I_{\ell}(\frg)$ denote the kernel of $\beta^\C$. Since $\mathcal{R}_{\ell}^\C(\frg)$ is
an algebra without nilpotents,  $I_{\ell}(\frg)$ is a radical
ideal.

{\it The main aim of this note is to conjecturally describe this ideal  $I_{\ell}(\frg)$.}
Before we describe our result and conjecture,
we  briefly describe the known results in this direction. Identify the
complexified representation  ring
$\mathcal{R}^\C(\frg)$
with the polynomial ring $\mathbb{C}[\chi_1,\ldots,\chi_r]$,
where $r$ is the rank of $\mathfrak{g}$ and $\chi_i$ denotes the
character of $V({\omega_i})$, the $i^{th}$ fundamental representation
of $\frg$. It is generally believed (initiated by the physicists) that there
exists an explicit potential function $F=F_\ell(\frg)$ (depending upon $\frg$ and $\ell$) in
$\mathbb{Z}[\chi_1,\ldots,\chi_r]$, coming from
representation theory of $\mathfrak{g}$,  with the property that the
ideal generated over the integers by the gradient of $F$, i.e.,  $\langle
\partial F/\partial \chi_1,\ldots, \partial
F/\partial \chi_r \rangle$, is precisely the kernel of $\beta$. This,
in particular, would imply  that the spec of $\mathcal{R}_{\ell}(\frg)$
is a complete intersection over $\mathbb{Z}$ inside the spec of the representation ring
$\mathbb{Z}[\chi_1,\ldots,\chi_r]$.  It may be remarked that, since  $I_{\ell}(\frg)$ is a radical
ideal of $\mathbb{C}[\chi_1,\ldots,\chi_r]$ with finite dimensional quotient,
 there exists an abstract potential function $F$ such that the
ideal generated by the gradient of $F$ over the complex numbers coincides with
 $I_{\ell}(\frg)$. However, an explicit construction of such an $F$ attuned to
 representation theory is only known in the case of $\frg=sl_{r+1}$ and $sp_{2r}$.

We recall the following result (cf.  \cite{G1}, \cite{G2},
\cite{BMRS} and \cite{BR}) obtained by first constructing explicitly
a potential function. Let $\chi_{\lambda}$ denote the character of
the irreducible representation $V(\lambda)$ of $\frg$.

\begin{theorem} (a) For $\frg={sl}_{r+1}$,
\[I_{\ell}({sl}_{r+1})=\langle \chi_{(\ell+1)\omega_1},
\chi_{(\ell+2)\omega_1},\dots, \chi_{(\ell+r)\omega_1}\rangle.\]

(b) For $\frg=sp_{2r}$,
\[I_{\ell}({sp}_{2r})=\langle
\chi_{\ell\omega_1+\omega_i};\;1\leq i \leq r
\rangle.\]

In particular, the ideals on the right side of both of (a) and (b) are radical.
\end{theorem}

For the remaining simple Lie algebras, we have the following theorem and the conjecture.

For any ideal $I$ in a ring $R$, we denote its radical ideal by
$\sqrt{I}$. In the following theorem, for any $\frg$, we take any
fundamental weight  $\omega_d$ with minimum Dynkin index. (Recall
from Section \ref{dynkin} that $\omega_d$ is unique  up to a diagram
automorphism except for $\frg$ of type $B_3$.) The list of such
$\omega_d$'s as well as the dual Coxeter number $\check{h}(\frg)$ of
any $\frg$ is given in Table 1 (Section 2).

\begin{theorem}
 (a) For $\frg$ of type $B_r (r\geq 3)$, $D_r (r\geq 4)$, $E_6$ or $E_7$,
\[
I_{\ell}(\mathfrak{g})\supseteq \sqrt{\langle \chi_{(\ell+1)\omega_d},
\chi_{(\ell+2)\omega_d},\ldots, \chi_{(\ell+
\check{h}(\frg)-1)\omega_d}\rangle}.\]
(For $\frg$ of type $B_3$, we must take $\omega_d=\omega_1$.)

 (b) For $\frg$ of type $G_2$,
\[I_{\ell}(G_2)\supseteq
\begin{cases}
\sqrt{\langle \chi_{(\ell+1)\omega_1}, \chi_{(\ell+2)\omega_1},
\chi_{\left(\frac{\ell+1}{2}\right)\omega_2}\rangle},& \text{if $\ell$ is odd,}\\
\sqrt{\langle \chi_{(\ell+1)\omega_1}, \chi_{(\ell+2)\omega_1},
\chi_{\omega_1+\frac{\ell}{2}\omega_2}\rangle},& \text{if $\ell$ is even.}\\
\end{cases}\]

(c)   For $\frg$ of type $F_4$,
\[I_{\ell}(F_4)\supseteq
\sqrt{\langle \chi_{(\ell+1)\omega_4},
\chi_{(\ell+2)\omega_4},\ldots, \chi_{(\ell+ 6)\omega_4}\rangle}\]

(d)  For $\frg$ of type  $E_8$,
\[I_{\ell}(E_8)\supseteq
\begin{cases} \sqrt{\langle \chi_{(\ell+2)\omega_8},
\chi_{(\ell+3)\omega_8},\ldots, \chi_{(\ell+29)\omega_8}\rangle},& \text{if $\ell$ is even,}\\
\sqrt{\langle \chi_{(\ell+2)\omega_8},
\chi_{(\ell+3)\omega_8},\dots, \chi_{(\ell+29)\omega_8},
\chi_{\frac{\ell+1}{2}\omega_8} \rangle},& \text{if $\ell$ is
odd.}\end{cases}\]

\end{theorem}

\begin{conjecture}  All the inclusions in (a)-(b) of the above theorem are, in fact,
 equalities for $\frg$ of type $B_r (r\geq 3), D_r (r\geq 4)$ and $G_2$.

In addition, we conjecture the following.

(a) For $\frg$ of type  $B_r (r\geq 3)$,
\begin{multline*}
\sqrt{\langle \chi_{(\ell+1)\omega_1},
\chi_{(\ell+2)\omega_1},\ldots, \chi_{(\ell+
\check{h}(\frg)-1)\omega_1}\rangle}= \\
{\langle \chi_{(\ell+1)\omega_1},
\chi_{(\ell+2)\omega_1},\ldots, \chi_{(\ell+
\check{h}(\frg)-1)\omega_1}, \chi_{\ell\omega_1+\omega_r}\rangle}.
  \end{multline*}

(b) For $\frg$ of type  $D_r (r\geq 4)$,
 \begin{multline*}\sqrt{\langle \chi_{(\ell+1)\omega_d},
\chi_{(\ell+2)\omega_d},\ldots, \chi_{(\ell+
\check{h}(\frg)-1)\omega_d}\rangle }= \\
{\langle \chi_{(\ell+1)\omega_d}, \chi_{(\ell+2)\omega_d},\ldots,
\chi_{(\ell+ \check{h}(\frg)-1)\omega_d},
\chi_{\ell\omega_1+\omega_{r-1}},
\chi_{\ell\omega_1+\omega_r}\rangle} .
\end{multline*}

\end{conjecture}

It is not clear to us (even conjecturally) how to explicitly construct potential
functions in these cases.

\noindent
{\bf Acknowledgement.} The second author was supported by an NSF grant. We thank
Peter Bouwknegt for some helpful correspondences.

\section{Preliminaries and Notation}

This section is devoted to setting  up the notation and recalling
some basic facts about affine Kac-Moody Lie algebras.

Let $G$ be a connected, simply-connected, simple algebraic group
over $\C$. We fix a Borel subgroup $B$ of $G$ and a maximal torus
$T\subset B$. Let $\frh$,$\frb$ and $\frg$ denote the Lie algebras of
$T$, $B$ and $G$ respectively.

Let $R=R(\frh,\frg)\subset \frh^*$ be
the root system; there is the root space decomposition $\frg=\frh
\oplus (\oplus_{\alpha \in R}\,\frg_{\alpha})$, $\frg_\alpha$ being the
 root space corresponding to the root $\alpha$.
The choice of $\frb$ determines a set of simple roots $\Delta=\{\alpha_1, \dots, \alpha_r\}$ of $R$,
where $r$ is the rank of $G$.
 For each root $\alpha$, denote by $\alpha^\vee$ the
unique element of $[\frg_{\alpha},\frg_{-\alpha}]$ such that
$\alpha(\alpha^\vee)=2$; it is called the {\it coroot} associated to the
root $\alpha$. Let $\frh_{\mathbb{R}}$ denote the real span of the
elements $\{\alpha_i^\vee, \alpha_i\in \Delta \}.$

Let $\{\omega_i \}_{1\leq i\leq r}$ be the set of fundamental
weights, defined as the basis of $\frh^*$ dual to
$\{\alpha_i^\vee\}_{1\leq i\leq r}$. We define the weight lattice
$P=\{ \lambda \in \frh^* : \lambda(\alpha_i^\vee) \in \mathbb{Z}\;
\forall \alpha_i \in \Delta \}$, and denote the set of dominant integral
weights by
$P_+$, i.e., $P_+:=\{ \lambda \in \frh^* : \lambda(\alpha_i^\vee) \in
\mathbb{Z}_+\;
\forall \alpha_i \in \Delta \}$, where $\mathbb{Z}_+$ is the set of
nonnegative integers.  The latter set parametrizes the set of isomorphism
classes of all the finite dimensional irreducible representations of $\frg$.
For $\lambda\in P_+$, let $V(\lambda)$
be the associated  finite dimensional irreducible
representations of $\frg$ with highest weight $\lambda$.  Let $\rho$ denote the sum of fundamental
weights, and $\check{h}=\check{h}(\frg):=1+\rho({\theta}^\vee)$ the {\it dual Coxeter
number}, where $\theta$ is the highest root of
$\frg$.

Let $(\;|\;)$ denote the Killing form on $\frg$ normalized such that
$({\theta}^\vee|{\theta}^\vee)=2$.  We will use the same notation for the restricted form on
$\frh$, and the induced form on $\frh^*$.
Let $W:=N_G(T)/T$ be the Weyl group of $G$, where $N_G(T)$ is the normalizer of $T$ in $G$.
Let $w_o$ denote
the longest element in the Weyl group.  For any $\lambda \in P$, the
dual of $\lambda$, denoted by $\lambda^*$, is defined to be $-w_o\lambda$.
For $\lambda\in P_+, \lambda^*$ is again in $P_+$ and, moreover,
 $V(\lambda)^*\simeq V(\lambda^*).$

\subsect{\bf Dynkin index.}\label{dynkin}
We recall the following definition from \cite{D}, \S2.
\begin{definition} Let $\mathfrak{g}_1$ and $\mathfrak{g}_2$ be two
(finite dimensional)  simple
Lie algebras and $\varphi : \mathfrak{g}_1 \to \mathfrak{g}_2$ be a
Lie algebra  homomorphism.  There exists a unique number $m_\varphi
\in \mathbb C$, called the {\it Dynkin index} of the homomorphism
$\varphi$, satisfying
\[(\varphi (x)|\varphi (y))
 = m_\varphi (x|y), \text{ for all }x,y \in \mathfrak{g}_1,\]
where $(\;|\;)$ is the normalized Killing form on
$\mathfrak{g}_1$ (and $\mathfrak{g}_2$) as defined above.
 Then, as proved by Dynkin  \cite{D}, Theorem 2.2, the Dynkin index is always a
nonnegative integer.

For a Lie algebra $\mathfrak g_1$ as above and a finite dimensional
representation $V$ of $\mathfrak g_1$, by the {\it Dynkin index}
$m_V$ of $V$, we mean the   Dynkin index of the Lie algebra
homomorphism $\rho_V : \mathfrak{g}_1 \to sl(V)$, where  $sl(V)$ is
the Lie algebra of traceless endomorphisms of $V$.
\end{definition}

For any simple  Lie algebra $\frg$, there is a fundamental weight
$\omega_d$  such that $m_{V(\omega_d)}$ divides $m_V$ for every
finite dimensional representation $V$. Moreover, such a $\omega_d$
is unique up to diagram automorphisms (except for $\frg$ of type
$B_3$). The following table gives the list of all such $\omega_d$'s
together with their indices $m_{V(\omega_d)}$ and also the dual
Coxeter number $\check{h}(\frg)$. The list of  $\omega_d$ and
$m_{V(\omega_d)}$ can be obtained from \cite{D}, Table 5 (see also
\cite{BK}, Proposition 2.3). We follow the indexing convention as in
\cite{B}.

\[ \begin{array}{cccc}
\mbox{Type of G}                     & \omega_d              & m_{V(\omega_d)} & \check{h}(\frg)\\
A_r \;(r\geq 1)                      & \omega_1, \omega_r    & 1 & r+1 \\
C_r \;(r\geq 2)                      & \omega_1              & 1 & r+1\\
B_r \;(r\geq 3)                      & \omega_1              & 2 & 2r-1\\
D_r \;(r\geq 4)                      & \omega_1              & 2 & 2r-2\\
G_2                                  & \omega_1              & 2 & 4\\
F_4                                  & \omega_4              & 6 & 9\\
E_6                                  & \omega_1, \omega_6 & 6 & 12\\
E_7                                  & \omega_7              & 12 & 18\\
E_8                                  & \omega_8              & 60 & 30.\\
&&&\\
& \mbox{Table 1.} &  &
\end{array} \]
We remark that for $B_3$, $\omega_3$ also satisfies $m_{V(\omega_3)}=2$; for $D_4$,
$\omega_3$ and $\omega_4$ both have
$m_{V(\omega_3)}=m_{V(\omega_4)}=2$.

\subsect{\bf Affine Kac-Moody Lie algebras.}
Let $\tilde{\frg}:=\frg \otimes \mathbb{C}((z))\oplus
\mathbb{C}K$ denote the (untwisted) affine Lie algebra
associated to $\frg$  (where $\mathbb{C}((z))$
denotes the field of Laurent series in one variable $z$), with
the Lie bracket
\[ [x\otimes f,y\otimes g]=[x,y]\otimes fg+ (x|y) Res_{z=0}(gdf)\cdot K
\;\mbox{and}\;[\tilde{\frg},K]=0,\] for $x,y \in \frg$ and $f,g \in
\mathbb{C}((z))$.

 The Lie algebra  $\frak{g}$ sits as a Lie subalgebra of
 $\tilde{\frak{g}}$ as $\frak{g} \otimes z^0$. The Lie
 algebra $\tilde{\frak{g}}$ admits a distinguished
 `parabolic' subalgebra
 $$\tilde{\frak{p}}:= \frak{g} \otimes \Bbb C [[z]] \oplus
 \Bbb C K.$$
 We also define its `nil-radical' $\tilde{\frak{u}}$ (which
 is an ideal of $\tilde{\frak{p}}$) by
 $$\tilde{\frak{u}}:= \frak{g} \otimes z \Bbb C [[z]],$$
 and its `Levi component' (which is a Lie subalgebra of
 $\tilde{\frak{p}})$
 $$\tilde{\frak{p}}^o:= \frak{g} \otimes z^0 \oplus \Bbb C
 K.$$
 Clearly (as a vector space)
 $$\tilde{\frak{p}}= \tilde{\frak{u}} \oplus
 \tilde{\frak{p}}^o.$$

\begin{subsect}{\bf Irreducible representations of
 $\tilde{\mathfrak g}$.}\label{2.5}
 Fix an irreducible  (finite dimensional) representation $V=V(\lambda)$
   of $\frak{g}$ and a positive integer  $\ell$.
 We define the associated {\it generalized Verma
 module} for $\tilde{\frak{g}}$ as
 $$M(V,  \ell)= U (\tilde{\frak{g}})
 \otimes_{U(\tilde{\frak{p}})} Y_{\ell} (V),$$
 where the $\tilde{\frak{p}}$-module $Y_{\ell} (V)$ has the underlying
 vector space same as  $V$ on which $\tilde{\frak{u}}$
 acts trivially, the central element $K$ acts via the scalar $\ell$
 and the action of $\frak{g}=\frak{g} \otimes z^0$ is
 via the $\frak{g}$-module structure on $V$.
Then,  $M(V, \ell)$ has a unique irreducible
 quotient denoted $L(V, \ell)$. Observe that $K$ acts by the constant
$\ell$ on  $M(V, \ell)$ and hence also on  $L(V, \ell)$. The constant by which $K$ acts on the
representation is called the {\it central charge} of the representation. Thus,  $L(V, \ell)$ has
central charge $\ell$.
  \end{subsect}

 \begin{definition}
  Consider the Lie subalgebra
 $\frak{r}^o$ of $\tilde{\frak{g}}$ spanned by
 $\{X_{-\theta} \otimes z, \theta^{\vee} \otimes 1, X_{\theta}
 \otimes z^{-1} \}$, where $X_{-\theta} $ (resp. $X_{\theta})$ is a
 non-zero root vector of $\frak{g}$ corresponding to the root
 $-\theta$ (resp. $\theta)$. Then the Lie algebra $\frak r^o$ is isomorphic with $sl(2)$.

 A $\tilde{\frak{g}}$-module $L$ is said to be {\it
 integrable} if every vector $v \in L$ is contained in a finite
 dimensional ${\frak g}$-submodule of $L$ and
 also $v$ is contained in a finite dimensional
 $\frak{r}^o$-submodule of $L$.

 Then, it follows easily from the  $sl
 (2)$-theory
 that the irreducible module $L(V(\lambda), \ell)$
  is integrable if and only if $\lambda\in P_\ell
:=\{\lambda \in P_+ : \lambda({\theta}^\vee)\leq
\ell\}$.
   \end{definition}

\section{An Introduction to the Fusion Ring}

Let $\mathcal{R}(\frg)$ denote the Grothendieck ring of finite
dimensional representations of $\frg$.  It is a free
$\mathbb{Z}$-module with basis $\{V({\lambda}):\;\lambda \in P_+\}$,
with the product structure induced from the tensor product of representations
 as follows:
\[V({\lambda})\otimes V({\mu})=\ds \sum_{\nu \in P_+} m_{\lambda ,
\mu}^{\nu} V({\nu}),\] with $m_{\lambda , \mu}^{\nu}=\dim
\Hom_{\frg}(V_{\lambda}\otimes V_{\mu}\otimes V_{\nu^*},
\mathbb{C}).$

The following is the most standard definition of the fusion ring.

  \begin{subsect}{\bf First Definition of the Fusion Ring.}
Fix a positive integer
$\ell$.  The {\it fusion ring} $\mathcal{R}_{\ell}(\mathfrak{g})$ (associated to $\mathfrak{g}$ and
the positive integer $\ell$)
 is a free $\mathbb{Z}$-module
with basis $\{ V({\lambda}): \lambda \in P_\ell \}$.
 The product structure, called the {\it fusion product}, is
defined as follows:

$$V({\lambda})\otimes ^F V{(\mu)}:= \ds \bigoplus_{\nu \in P_{\ell}}
n_{\lambda,\mu}^{\nu}(\ell) V(\nu),$$ where $n_{\lambda,\mu}^{\nu}(\ell)$ is
the dimension of the space $V_{\mathbb{P}^1}^\dag(\lambda,\mu,\nu^*)$ of {\it conformal blocks} on
$\mathbb{P}^1$ with three distinct marked points $\{x_0,x_1,x_2\}$ and the weights
$\lambda,\mu,\nu^*$ attached to them with central charge $\ell$.
Recall that
$$V_{\mathbb{P}^1}^\dag(\lambda,\mu,\nu^*):=
\Hom_{\frg\otimes\mathcal{O}[U]}\bigl(L(V({\lambda}),\ell)\otimes L(V({\mu}),\ell)\otimes
L(V({\nu^*}),\ell),\mathbb{C}\bigr),$$ where $U:=\mathbb{P}^1\setminus\{x_0,x_1,x_2\}$ and
$\mathcal{O}[U]$ denotes
the ring of $\Bbb C$-valued regular functions on the affine curve $U$.
The
action of $\frg\otimes\mathcal{O}[U]$
on $L(V({\lambda}),\ell)\otimes L(V({\mu}),\ell)\otimes
L(V({\nu^*}),\ell)$ is given by
\begin{multline*}(x\otimes f)(v_1\otimes v_2\otimes v_3):=\\
((x\otimes(f)_{x_0})v_1)\otimes v_2\otimes
v_3 + v_1\otimes ((x\otimes(f)_{x_1})v_2)\otimes v_3
 +v_1\otimes v_2 \otimes ((x\otimes(f)_{x_2})v_3),
\end{multline*}
 for $x\in \frg$ and $f\in \mathcal{O}[U]$, where $(f)_{x_i} \in \Bbb C((z))$ denotes
 the Laurent series expansion of $f$ at $x_i$. The numbers
$n_{\lambda,\mu}^{\nu}(\ell)$ are  given by the celebrated Verlinde
formula \cite{V}.

It is clear from its definition that the product $\otimes^F$ is
commutative.  It is also associative as a result of the
`factorization rules' of \cite{TUY}.

We also remark that the canonical map
\[p:V_{\mathbb{P}^1}^\dag(\lambda,\nu,\mu^*)\to \Hom_{\frg}(V(\lambda)\otimes V(\mu)\otimes
V(\nu^*), \mathbb{C})\] induced from the natural inclusion
$V(\lambda)\otimes V(\mu)\otimes V(\nu^*)\hookrightarrow
L(V({\lambda}),\ell)\otimes L(V({\mu}),\ell)\otimes
L(V({\nu^*}),\ell)$ is an injection
\cite{SU}.  In particular,  the
inequality $n_{\lambda,\mu}^{\nu}(\ell)\leq m_{\lambda,\mu}^{\nu}$ holds.
  \end{subsect}

\begin{subsect}{\bf Second definition.}
 Let $\tilde{\mathcal G}$ be the affine Kac-Moody group
 associated to the Lie algebra $\tilde{\frak{g}}$ and
 $\tilde{\mathcal P}$ its parabolic subgroup (corresponding
 to the Lie subalgebra $\tilde{\frak{p}}$) (cf. \cite{Ku1}, Chapter 13). Then,
 $X=\tilde{\mathcal G}/\tilde{\mathcal P}$ is a projective  ind-variety.
 Now, given a finite dimensional algebraic representation $V$ of
 $\tilde{\mathcal P}$, we can consider the associated
 homogeneous vector bundle $\mathcal V$ on $X$ and the
 corresponding Euler-Poincar\'e characteristic (which
 is a virtual $\tilde{\mathcal G}$-module)
 $$
 \mathcal X(X,\mathcal V ) :=  \sum_i (-1)^i H^i (X,\mathcal V ).
 $$
Recall that $H^i(X, \mathcal V)$ is determined in \cite{Ku1}, Chapter 8.

 For any positive integer $\ell $ and
 $\lambda,\mu \in P_\ell $, define
 $$
 [L(V(\lambda),\ell) \otimes^{\bullet}
 L(V(\mu),\ell ) ]^\ast
 \cong
 \mathcal X (X,\mathcal V ),
 $$
 as virtual $\tilde{\mathcal G}$-modules,
 where the $\tilde{\mathcal P}$-module  $V := (Y_\ell(V(\lambda) \otimes V(\mu)))^\ast $
 (cf. \S\ref{2.5} for the notation $Y_\ell$). Writing
\[\mathcal X (X,\mathcal V )\cong \oplus_{\nu \in P_\ell}\,d^\nu_{\lambda, \mu}\,L(V(\nu),\ell )^* ,\]
we get another definition of the fusion product as follows:
$$V(\lambda)\otimes^{\bullet} V{(\mu)}:= \ds \bigoplus_{\nu \in P_{\ell}}
d_{\lambda,\mu}^{\nu} V(\nu).$$

By a result of Faltings \cite{F}, Appendix (see also Kumar \cite{Ku}, Theorem 4.2), the two products $\otimes^F$ and
$\otimes^{\bullet}$
coincide for $G$ of type $A_r,B_r,C_r,D_r$ and $G_2$. In fact, these
 products coincide
for any $G$ as proved by Teleman (cf.
\cite{T}, Theorem 0 together with  \cite{Ku}, Lemma 4.1).
  \end{subsect}

\begin{subsect}{\bf Third definition.} Consider the Lie subalgebra $\mathfrak s$ of
 $\frg$ spanned by
$X_\theta,X_{-\theta}$ and $\theta^\vee$, where $X_{\pm \theta}$ are the root vectors
corresponding to the roots $\pm \theta$ respectively, choosen so that
 $(X_\theta | X_{-\theta})=1$. Then, $\mathfrak s$ is isomorphic with $sl_2$
under the isomorphism $\phi: sl_2 \to \mathfrak s$ taking $X\mapsto X_\theta,
Y\mapsto X_{-\theta}, H\mapsto \theta^\vee$, where $X,Y,H$ is the standard
 basis of $sl_2$. For any finite dimensional $\mathfrak s$-module $V$
 and nonnegative integer $d$, let $V^{(d)}$ denote the isotypical component of
 the $\mathfrak s$-module $V$ corresponding to the highest weight $d$
(where we follow the convention that the highest weight of an irreducible
 $sl_2$-module is one less than its dimension).
  \end{subsect}
The following proposition can be found in  \cite{Bea}, Proposition 7.2.

\begin{proposition} For any $\lambda,\mu \in P_\ell$, the fusion product
 $V(\lambda)\otimes ^F V{(\mu)}$ is the isomorphism class of the ($\frg$-module) quotient
of $V{(\lambda)}\otimes  V{(\mu)}$ by the $\frg$-submodule generated
by the
 $\mathfrak s$-module
$$\bigoplus_{p+q+d >2\ell} \,\bigl(V(\lambda)^{(p)}\otimes
 V(\mu)^{(q)}\bigr)^{(d)}.$$
\end{proposition}

\begin{subsect}{\bf The homomorphism $\beta : \mathcal{R}(\frg) \to \mathcal{R}_{\ell}(\frg)$.}\label{3.5}
We first introduce the affine Weyl group $W_{\ell}$. As in Section 2,
let $W$ be the (finite) Weyl group of $G$ which acts naturally on $P$ and hence
 also on $P_{\Bbb R}:=P\otimes_{\Bbb Z}\Bbb R$. Let $W_\ell$ be the group of
 affine transformations of $P_{\Bbb R}$ generated by $W$ and the translation
 $\lambda\mapsto \lambda+(\ell+\check{h})\theta$. Then, $W_\ell$ is the semi-direct product
 of $W$ by the lattice $(\ell+\check{h})Q^{\text{long}}$, where
 $Q^{\text{long}}$ is the sublattice of $P$ generated by the long roots. For any root
 $\alpha\in R$ and $n\in \Bbb Z$, define the {\it affine wall}
 \[H_{\alpha,n}=\{\lambda\in P_{\Bbb R}: (\lambda|\alpha)=n(\ell+\check{h})\}.\]
 The closures of the connected components of $P_{\Bbb R}\setminus
 \bigl(\cup_{\alpha\in R, n\in \Bbb Z}\,H_{\alpha,n}\bigr)$ are called the {\it alcoves}.
 Then, any alcove is a fundamental domain for the action of $W_\ell$. The
 {\it fundamental alcove} is by definition
 \[A^o=\{\lambda\in P_{\Bbb R}: \lambda(\alpha_i^\vee)\geq 0\, \forall\,
 \alpha_i \in \Delta, \,\text{and}\, \lambda(\theta^\vee)\leq \ell+\check{h}\}.\]
  \end{subsect}

 \begin{definition}\label{kerpi}  Define the $\Bbb Z$-linear map
$\beta :  \mathcal{R}(\frg) \to   \mathcal{R}_{\ell}(\frg)$ as follows.  Let $\lambda \in P_+$.
  If $\lambda+\rho$ lies
on an affine wall, then $\beta(V({\lambda}))=0$.  Otherwise, there is a unique
$\mu\in P_\ell$ and $w\in W_\ell$ such that  $\lambda+\rho=w(\mu+\rho)$. In this case,
define
$\beta(V({\lambda}))=\varepsilon(w)V(\mu)$, where $\varepsilon(w)$ is the sign of
the affine Weyl group element $w$.
\end{definition}

The following theorem follows from the equivalence of the first two
definitions of the fusion product given above together with \cite{Ku}, Lemma 3.3
and the Kac-Moody analogue of the Borel-Weil-Bott theorem (cf. \cite{Ku1},
 Corollary 8.3.12).

\begin{theorem}  The $\Bbb Z$-linear map
$\beta$ defined above is, in fact, an algebra homomorphism with respect to
the fusion product on $ R_\ell(G)$.
\end{theorem}

\section{A Conjectural Presentation of the Fusion Algebra}

We  consider the complexification $\mathcal{R}^{\mathbb{C}}(\frg)$ of $\mathcal{R}(\frg)$,
 and similarly consider the  complexification
 $\mathcal{R}^{\mathbb{C}}_\ell(\frg)$ of $\mathcal{R}_{\ell}(\frg)$. We will
 denote the character of the $i^{th}$ fundamental
representation of $\frg$ by $\chi_i$ (instead of $\chi_{\omega_i}$).
As is well known, $\mathcal{R}(\frg)$ is identified
with the polynomial ring $\mathbb{Z}[\chi_1,\ldots,\chi_r]$. Moreover,
the $\Bbb C$-algebra  $\mathcal{R}^{\mathbb{C}}(\frg)$ is identified with
the affine coordinate ring $\Bbb C[T/W]$ of the quotient $T/W$. This
identification
is obtained by taking the character of any virtual representation in
 $\mathcal{R}(\frg)$. Similarly, the finite dimensional algebra
 $\mathcal{R}^{\mathbb{C}}_\ell(\frg)$ is identified with the affine
coordinate ring
$\Bbb C[T^{\text{reg}}_{\ell}/W]$ of the reduced subscheme
$T^{\text{reg}}_{\ell}/W$ of $T/W$ by taking the character values on the
finite set $T^{\text{reg}}_{\ell}/W$ (cf. \cite{Bea}, \S9), where
\[T_\ell:=\{t\in T: e^\alpha(t)=1 \, \forall \alpha\in (\ell+
\check{h})Q^{\text{long}}\},\] and $T^{\text{reg}}_{\ell}$ is the
subset of regular elements of $T_\ell$. Moreover, the homomorphism
$\beta^{\Bbb C} : \mathcal{R}^{\mathbb{C}}(\frg)\to
\mathcal{R}^{\mathbb{C}}_\ell(\frg)$
 obtained from the
 complexification of $\beta$ defined in Definition \ref{kerpi} (under the above identifications)
 corresponds to the inclusion of
 $T^{\text{reg}}_{\ell}/W \hookrightarrow T/W$.
Let   $I_{\ell}(\frg)$ denote the kernel of $\beta^{\Bbb C}$.
  Since
 $T^{\text{reg}}_{\ell}/W$ is reduced,  $I_{\ell}(\frg)$  is a
radical ideal. Identifying  $\mathcal{R}^{\Bbb C}(\frg)$
with  $\mathbb{C}[\chi_1,\ldots,\chi_r]$ as above, we can think of
 $I_{\ell}(\frg)$ as an ideal inside  $\mathbb{C}[\chi_1,\ldots,
\chi_r]$. Thus,
\[ \mathcal{R}^{\mathbb{C}}_\ell(\frg)\simeq \mathbb{C}[\chi_1,\ldots,
\chi_r]/I_{\ell}(\frg).\]

We begin by recalling the following result on the presentation of fusion algebras.
The (a)-part of the result is due to Gepner \cite{G1} and for the (b)-part, see
\cite{G2} and \cite{BMRS}. In fact, both the parts are obtained by first constructing
explicitly a potential function. Moreover, the following presentations are valid even over
$\Bbb Z$ (cf. \cite{BR}).

\begin{theorem} (a) For $\frg={sl}_{r+1}$,
\[I_{\ell}({sl}_{r+1})=\langle \chi_{(\ell+1)\omega_1},
\chi_{(\ell+2)\omega_1},\ldots, \chi_{(\ell+r)\omega_1}\rangle.\]
(Recall that for ${sl}_{r+1}, \check{h}=r+1.$)

(b) For $\frg=sp_{2r}$,
\[I_{\ell}({sp}_{2r})=\langle
\chi_{\ell\omega_1+\omega_i};\;1\leq i \leq r
\rangle.\]

In particular, the ideals on the right side of both of (a) and (b) are radical.
\end{theorem}

For any ideal $I$ in a ring
 $R$, we denote its radical ideal by  $\sqrt{I}$. Recall that the complete list
 of the fundamental weights $\omega_d$ such that the Dynkin index $m_{V(\omega_d)}$ is
 minimum (as well as the dual Coxeter numbers $\check{h}$ of any $\frg$) are given in Table 1.
 In the following theorem, for any $\frg$, we take any one $\omega_d$ with minimum Dynkin index.
(Recall from Section \ref{dynkin} that $\omega_d$ is unique  up to a
diagram automorphism except
 for $\frg$ of type $B_3$.)

\begin{theorem} \label{main} Let $\ell$ be any positive integer.

 (a) For $\frg$ of type $B_r (r\geq 3)$, $D_r (r\geq 4)$, $E_6$ or $E_7$,
\[
I_{\ell}(\mathfrak{g})\supseteq \sqrt{\langle \chi_{(\ell+1)\omega_d},
\chi_{(\ell+2)\omega_d},\ldots, \chi_{(\ell+
\check{h}(\frg)-1)\omega_d}\rangle}.\]
(For $\frg$ of type $B_3$, we must take $\omega_d=\omega_1$.)

 (b) For $\frg$ of type $G_2$,
\[I_{\ell}(G_2)\supseteq
\begin{cases}
\sqrt{\langle \chi_{(\ell+1)\omega_1}, \chi_{(\ell+2)\omega_1},
\chi_{\left(\frac{\ell+1}{2}\right)\omega_2}\rangle},& \text{if $\ell$ is odd,}\\
\sqrt{\langle \chi_{(\ell+1)\omega_1}, \chi_{(\ell+2)\omega_1},
\chi_{\omega_1+\frac{\ell}{2}\omega_2}\rangle},& \text{if $\ell$ is even.}\\
\end{cases}\]

(c)   For $\frg$ of type $F_4$,
\[I_{\ell}(F_4)\supseteq
\sqrt{\langle \chi_{(\ell+1)\omega_4},
\chi_{(\ell+2)\omega_4},\ldots, \chi_{(\ell+ 6)\omega_4}\rangle}.\]

(d)  For $\frg$ of type  $E_8$,
\[I_{\ell}(E_8)\supseteq
\begin{cases} \sqrt{\langle \chi_{(\ell+2)\omega_8},
\chi_{(\ell+3)\omega_8},\ldots, \chi_{(\ell+29)\omega_8}\rangle},& \text{if $\ell$ is even,}\\
\sqrt{\langle \chi_{(\ell+2)\omega_8},
\chi_{(\ell+3)\omega_8},\dots, \chi_{(\ell+29)\omega_8},
\chi_{\frac{\ell+1}{2}\omega_8} \rangle},& \text{if $\ell$ is
odd.}\end{cases}\]

\end{theorem}

\begin{proof}
Using the tables in Bourbaki \cite{B}, in the cases (a), it is easy
to verify that for $\lambda=(\ell+m)\omega_d$, $1\leq m\leq
\check{h}(\frg)-1$, there exists an affine wall $H_{\mu, n}$
determined by a positive root $\mu$ of $\frg$ such that
$\lambda+\rho$ lies on it. Then, by the definition of $\beta$,
 $\beta(V({\lambda}))$ is trivial.

We give below some of the details.
\medskip

 \underline{$B_r$}: In this case, we have $\omega_d=\omega_1$ and
 $\check{h}=2r-1$. Observe that $((\ell+m)\omega_1+\rho|\mu)=\ell+2r-1$,
where
\[\begin{cases}
\mu=\sum_{1\leq k\leq m}\alpha_k +2\sum_{1+m\leq k\leq r} \alpha_k
&{\rm if}\;\;1\leq m \leq r-1, \\
\mu=\sum_{1\leq k<2r-m}\alpha_k &{\rm if} \;\; r\leq m \leq 2r-2.
\end{cases}\]
Thus, $(\ell+m)\omega_1+\rho$ lies on the affine wall
determined by $\mu$ as above.

 Similarly,  $\ell
\omega_1+\omega_r+\rho$ lies on the wall determined by the highest
root $\theta$ of $B_r$ (see the next conjecture, part (a)).

\bigskip
\underline{$D_r$}: In this case, we take $\omega_d=\omega_1$ and we have
 $\check{h}=2r-2$.  Observe that
$((\ell+m)\omega_1+\rho|\mu)=\ell+2r-2$,  where
\[\begin{cases}
\mu=\sum_{1\leq k\leq m}\alpha_k+ (2\sum_{m+1\leq k<r-1}\alpha_k) +\alpha_{r-1}+\alpha_r  &{\rm if} \;\; 1\leq m \leq r-2,\\
\mu=\sum_{1\leq k<2r-1-m}\alpha_k &{\rm if}\;\; r-1\leq m \leq 2r-3.
\end{cases}\]
Since any $\omega_d$ is obtained from $\omega_1$ by a diagram automorphism
in this case, the result
 holds for any $\omega_d$.

A similar calculation shows that both $\ell \omega_1+\omega_{r-1}+\rho$
and $\ell \omega_1+\omega_r+\rho$ lie on the affine wall determined
by the highest root $\theta$ of $D_r$ (see the next conjecture, part (b)).

\bigskip

\underline{$G_2$}: In this case, we have  $\omega_d=\omega_1$ and
 $\check{h}=4$. Observe that $(\ell+1)\omega_1+\rho$ (resp.~
$(\ell+2)\omega_1+\rho$) lies on the affine wall determined by
$\theta=3\alpha_1+2\alpha_2$ (resp.~$3\alpha_1+\alpha_2$).  Moreover,
$\frac{\ell+1}{2}\omega_2+\rho$ (for odd $\ell$) and
$\omega_1+\frac{\ell}{2}\omega_2+\rho$ (for even $\ell$) lie on the affine wall
determined by the root $\theta$.

\underline{$F_4$}: In this case, we have  $\omega_d=\omega_4$ and
 $\check{h}=9$. Observe that $(\ell+1)\omega_4+\rho,\dots,  (\ell+6)\omega_4+\rho$
 lie on the affine wall determined by the roots $\theta=2\alpha_1+3\alpha_2+4\alpha_3+2\alpha_4,
 \, \theta-\alpha_1,\, \theta-\alpha_1
 -\alpha_2,\, \theta-\alpha_1-\alpha_2-2\alpha_3,\, \theta-\alpha_1-2\alpha_2-2\alpha_3,
 \,\theta-2\alpha_1-2\alpha_2-2\alpha_3$ respectively.

The corresponding calculation  for the $E$ series is more involved
and is left to the reader. This proves the theorem.
\end{proof}

We would like to make the following conjecture.
\medskip

\begin{conjecture}
 All the inclusions in (a)-(b) of the above theorem are, in fact,
 equalities for $\frg$ of type $B_r (r\geq 3), D_r (r\geq 4)$ and $G_2$.

 In addition, we conjecture the following:

(a) For $\frg$ of type  $B_r (r\geq 3)$,
  \begin{multline*}
\sqrt{\langle \chi_{(\ell+1)\omega_1},
\chi_{(\ell+2)\omega_1},\ldots, \chi_{(\ell+
\check{h}(\frg)-1)\omega_1}\rangle} = \\
{\langle \chi_{(\ell+1)\omega_1},
\chi_{(\ell+2)\omega_1},\ldots, \chi_{(\ell+
\check{h}(\frg)-1)\omega_1},
 \chi_{\ell\omega_1+\omega_r}\rangle}.
  \end{multline*}

(b) For $\frg$ of type  $D_r (r\geq 4)$:
 \begin{multline*}
\sqrt{\langle \chi_{(\ell+1)\omega_d},
\chi_{(\ell+2)\omega_d},\ldots, \chi_{(\ell+
\check{h}(\frg)-1)\omega_d}\rangle}=\\
{\langle \chi_{(\ell+1)\omega_d},
\chi_{(\ell+2)\omega_d},\ldots, \chi_{(\ell+
\check{h}(\frg)-1)\omega_d},
 \chi_{\ell\omega_1+\omega_{r-1}},
\chi_{\ell\omega_1+\omega_{r}}\rangle}.\end{multline*}
\end{conjecture}

\noindent {\bf 4.4. Question.} One may ask if the inclusions in the
above theorem \ref{main} for $\frg$ of type $F_4, E_6, E_7$ and
$E_8$ are equalities as well.

\medskip
\medskip
Address: Department of Mathematics, UNC at Chapel Hill, Chapel Hill, NC 27599-3250, USA
(email:  boysal@math.jussieu.fr and shrawan@email.unc.edu)


\begin{thebibliography}{999}

\bibitem{B} N.~Bourbaki. Groupes et Alg\`ebres de Lie, Chap. 4--6. Masson, Paris (1981).

\bibitem{Bea}  A.~Beauville. Conformal blocks, fusion rules and the Verlinde formula.
Proceedings of the Hirzebruch 65 Conference on Algebraic Geometry
(1996) 75--96.

\bibitem{BMRS} M.~Bourdeau, E.~Mlawer, H.~Riggs and H.~Schnitzer.
Topological Landau-Ginzburg matter from $Sp(N)_k$ fusion rings.
Modern Physics Letters A {\bf 7} (1992)  689--700.

\bibitem{BR} P.~Bouwknegt and D.~Ridout.  Presentations of Wess-Zumino-Witten fusion rings.
Rev.~Math.~Phys. {\bf 18} (2006)  201--232.


\bibitem{BK}  A.~Boysal and S.~Kumar. Explicit determination of the Picard group of
moduli spaces of semistable $G$-bundles on curves.  Math. Annalen {\bf 332} (2005)
823--842.

\bibitem{D} E.B.~Dynkin, Semisimple subalgebras of semisimple Lie algebras. Am. Math. Soc. Transl.
(Ser. II) {\bf 6} (1957) 111--244.

\bibitem{F}  G.~Faltings. A proof for the Verlinde formula.  J.~Alg.~Geom. {\bf 3} (1994)
347--374.

\bibitem{G1} D.~Gepner. Fusion rings and geometry.
Comm. Math. Phys.   {\bf 141}
 (1991) 381--411.

\bibitem{G2} D.~Gepner and A. Schwimmer. Symplectic fusion rings and their metric.
Nuclear Physics  B {\bf 380} (1992) 147--167.

\bibitem{K} V.~Kac. Infinite Dimensional Lie Algebras. Cambridge University Press (1992).

\bibitem{Ku} S.~Kumar.  Fusion product of positive level
representations.  Lecture notes in Pure and Applied Mathematics {\bf 184}
(1997) 253--259.

\bibitem{Ku1} S.~Kumar. Kac-Moody groups, their flag varieties and representation theory.
Progress in Mathematics vol. {\bf 204}, Birkhauser (2002).

\bibitem{T} C.~Teleman. Lie algebra cohomology and the fusion rules. Comm. Math. Phys.
{\bf 173} (1995) 265--311.

\bibitem{TUY} A.~Tsuchiya, K.~ Ueno and Y.~Yamada. Conformal field theory on universal
family of stable curves with gauge symmetries. Adv.~Stud.~Pure Math. {\bf 19}
(1989)  459--565.

\bibitem{SU} Y.~Shimizu and  K.~ Ueno. Advances in Moduli Theory.
Tranlations of Mathematical Monographs  {\bf 206}, AMS (2002).

\bibitem{V} E.~Verlinde. Fusion rules and modular transformations in 2d conformal field theory.
Nuclear Physics B {\bf 300} (1988) 360--376.

\end{thebibliography}
\end{document}